\documentclass[11pt]{amsart}
\usepackage{amscd,amssymb,amsmath,latexsym,enumerate}
\usepackage[mathscr]{euscript}
\usepackage{epsfig}
\usepackage{color}
\definecolor{gray}{rgb}{0.93,0.93,0.93}
\usepackage{tikz}
\usetikzlibrary{decorations,arrows}
\usetikzlibrary{decorations.pathmorphing}
\usepgflibrary{decorations.pathreplacing} 
\usepackage{pgf}
\usepackage{pgfkeys}
\usepackage{pgffor}
\usepackage{pgfpages}

\newtheorem{theo}{Theorem}
\newtheorem{defini}{Definition}
\newtheorem{proposi}{Proposition}
\newtheorem{lemma}{Lemma}

\newtheorem{rem}{Remark}
\newtheorem{exam}{Example}

\newcommand{\Ff}{{\mathcal F}}

\newcommand{\Hh}{{\mathcal H}}

\newcommand{\Kk}{{\mathcal K}}

\newcommand{\Uu}{{\mathcal U}}

\newcommand{\CM}{{\mathbb C}}

\newcommand{\NM}{{\mathbb N}}
\newcommand{\QM}{{\mathbb Q}}

\newcommand{\RM}{{\mathbb R}}

\newcommand{\TM}{{\mathbb T}}

\newcommand{\ZM}{{\mathbb Z}}

\newcommand{\gG}{{\mathfrak g}}

\newcommand{\hf}{{\widehat{f}}}

\newcommand{\hg}{{\widehat{g}}}

\newcommand{\hp}{{\widehat{p}}}

\newcommand{\cs}{{\mathscr C}}

\newcommand{\ps}{{\mathscr P}}

\newcommand{\ts}{{\mathscr T}}

\newcommand{\ws}{{\mathscr W}}

\newcommand{\Css}{$C^{\ast}$-algebras }       %%C*-algebras w. end space
\newcommand{\CsS}{$C^{\ast}$-algebras}             %%C*-algebras
\newcommand{\dist}{\mbox{\rm dist}}                %%distance
\newcommand{\Lip}{\mbox{\rm Lip}}	           %%Lipschitz constant
\newcommand{\dil}{\mbox{\rm dil}}	  	   %%dilation
\newcommand{\Hol}{\mbox{\rm Hol}}	  	   %%H\"older constant

\begin{document}

\title{Continuity of the spectrum of a field of self-adjoint operators}
\thanks{Work supported in part by NSF Grant DMS-1160962}
\author{Siegfried Beckus, Jean Bellissard}

\address{Mathematisches Institut\\
Friedrich-Schiller-Universit\"at Jena\\
07743, Jena, Germany}
\email{siegfried.beckus@uni-jena.de}

\address{Georgia Institute of Technology\\
School of Mathematics\\
Atlanta GA 30332-0160}
\email{jeanbel@math.gatech.edu}

\begin{abstract}
Given a family of self-adjoint operators $(A_t)_{t\in T}$ indexed by a parameter $t$ in some topological space $T$, necessary and sufficient conditions are given for the spectrum $\sigma(A_t)$ to be Vietoris continuous  with respect to $t$. Equivalently the boundaries and the gap edges are continuous in $t$. If $(T,d)$ is a complete metric space with metric $d$, these conditions are extended to guarantee H\"older continuity of the spectral boundaries and of the spectral gap edges. As a corollary, an upper bound is provided for the size of closing gaps. 
\end{abstract}

%%%%%%%%%%%%%%%%%%%%%%%%%%%%%%%%%%%%%%%%%%%%%%%%%%%%%%%%%%%%%%%%%%%%

\maketitle
%\tableofcontents

%%%%%%%%%%%%%%%%%%%%%%%%%%%%%%%%%%%%%%%%%%%%%%%%%%%%%%%%%%%%%%%%%%%%
\section{Introduction}
\label{ccs15.sect-Introduction}

\noindent Given a family of self-adjoint operators $(A_t)_{t\in T}$ indexed by a parameter $t$ in some topological space $T$, what conditions are needed for this family to insure that the spectrum $\sigma(A_t)$ varies continuously with $t$~? In this work we provide an answer using analytic methods.

\vspace{.1cm}

\noindent This problem is usually solved by using continuous fields of Banach spaces, in particular of Hilbert spaces and \CsS. This concept, initially proposed by Tomiyama \cite{TT61,To62}, was further developed by Dixmier and Douady \cite{DD63}. A suggested reference on this topics is the book by Dixmier \cite{Dix64}. Unlike fiber bundles, a continuous field of Banach spaces may have pairwise non isomorphic fibers. The field of irrational rotation algebras over the interval $[0,1]$ is a typical example \cite{Rie81}. Nevertheless one of the most powerful consequences of the continuity is that a continuous field of normal elements of a field of \Css admits a spectrum varying continuously with the parameter. In this work we seek to help clarify why this occurs. Moreover, we provide a proof of the result without the machinery of continuous fields of \CsS. 

%%%%%%%%%%%%%%%%%%%%%%%%%%%%%%%%%%%%%%%%%%%%%%%%%%%%%%%%%%%%%%%%%%%%
 \subsection{An Example: the Almost Mathieu Operator}
 \label{ccs15.ssect-AM}

\noindent To illustrate the main difficulty, let $H_t$, where $t \in [0,1]=T$, be the Almost Mathieu model acting on $\Hh=\ell^2(\ZM)$ as follows
%%%%%%%%%%%%%%%%
\begin{equation}
\label{ccs15.eq-am}
H_t\psi(n)=
 \psi(n+1)+ \psi(n-1)+ 2\mu \cos{2\pi(nt+ \theta)} \psi(n)\,,
  \hspace{2cm}
   \psi\in\ell^2(\ZM)\,,\;\;n\in\ZM.
\end{equation}
%%%%%%%%%%%%%%%%

\noindent In this definition $\theta$ is a fixed parameter and $\mu>0$. It is clear that $H_t$ is strongly continuous in $t\in T$. On the other hand, if $t\neq s$ are irrational and rationally independent, it is easy to check that $\|H_t-H_s\|=2\mu$, so that this family is not norm continuous. However, it has been shown \cite{AS85,Be94} that the norm of $\|p(H_t)\|$ is (Lipschitz) continuous in $t$ for each polynomial $p$ and it was deduced from this that the gap edges of the spectrum of $H_t$ were also (Lipschitz) continuous as long as the gap does not close. Near the points where a gap closes, the gap edges are only H\"older continuous of exponent $1/2$, c.f. \cite{RB90,HS90}.

%%%%%%%%%%%%%%%%%%%%%%%%%%%%%%%%%%%%%%%%%%%%%%%%%%%%%%%%%%%%%%%%%%%%
 \subsection{Main Results}
 \label{ccs15.ssect-main}

\noindent The general formulation of the problem requires the following assumptions:

\begin{itemize}
  \item[F1)] $T$ is a topological space.

  \item[F2)] $\Hh_t$ is a Hilbert space for each $t\in T$.

  \item[F3)] For each $t\in T$, $A_t$ is a linear, possibly unbounded, self-adjoint operator on $\Hh_t$.
\end{itemize}

\noindent The family $\Hh=(\Hh_t)_{t\in T}$ meeting the assumption [F2] is called a field of Hilbert spaces. Similarly, a family $A=(A_t)_{t\in T}$ satisfying the assumption [F3] is called a field of self-adjoint operators.

%%%%%%%%%%%%%%%%
\begin{defini}
\label{ccs15.def-p2}
Under the assumptions [F1-F3] above, a field $A$ of self-adjoint, bounded operators is called p2-continuous if the maps $\Phi_p: t\in T\mapsto \|p(A_t)\|$ are continuous whenever $p$ is a polynomial in one variable with real coefficients and degree at most $2$. 
\end{defini}
%%%%%%%%%%%%%%%%

\noindent The spectrum of a normal closed operator is given by a closed subset of the complex plane. If the operator is self-adjoint its spectrum is a closed subset of the real line. Consequently, in order to describe the continuity of the spectrum, one needs a topology on the set $\cs(X)$ of closed subsets of a given topological Hausdorff space $X$. Such topologies are also known under the names of {\em hyperspace} or {\em hit-and-miss topology} \cite{Beer,LP94}. It seems that the first example of such a topology was provided by the Hausdorff metric in \cite{Ha14} (see English versions in \cite{Ha35,Ha62}; see also \cite{Ku68,Mu75,CV77,Ba88}), whenever $X$ is a complete metric space. In 1922, Vietoris \cite{Vi22} (see also \cite{CV77}) proposed such a topology for any topological space $X$. He showed that if $X$ is a complete metric space, then his topology coincides with the one defined by the Hausdorff metric. In 1962, Fell \cite{Fe62} defined a slightly different topology for which $\cs(X)$ is always compact but not always Hausdorff, unless $X$ is locally compact. The Fell and the Vietoris topologies coincide if $X$ is compact and Hausdorff. 

\vspace{.1cm}

\noindent The first main result is the following.

%%%%%%%%%%%%%%%%
\begin{theo}
\label{ccs15.th-main1}
Let $A=(A_t)_{t\in T}$ be a family of self-adjoint, bounded operators satisfying the assumptions [F1-F3]. Then the following are equivalent: 

(i)  the spectrum $\sigma(A_t)$ is a Vietoris continuous function of $t$, 

(ii) the spectral edges are continuous, 

(iii) the field $A$ is p2-continuous.
\end{theo}
%%%%%%%%%%%%%%%%

\noindent The precise concept of spectral edges is described in Section~\ref{ccs15.ssect-gapedge}. The previous result is valid only for fields of bounded self-adjoint operators. What about unbounded ones?

%%%%%%%%%%%%%%%%
\begin{defini}
\label{ccs15.def-res}
Let $A=(A_t)_{t\in T}$ be a field of (not necessarily bounded) self-adjoint operators over $T$. Such a field will be called R-continuous if the norm of its resolvent is continuous. Precisely,  this means that for every $z\in \CM\setminus \RM$, the map $t\in T\mapsto \|(z-A_t)^{-1}\|\in [0,\infty)$ is continuous.
\end{defini}
%%%%%%%%%%%%%%%%

%%%%%%%%%%%%%%%%
\begin{theo}
\label{ccs15.th-mainR}
Let $A=(A_t)_{t\in T}$ be a field of self-adjoint operators. Then the following are equivalent

(i) $A$ is R-continuous,

(ii) the spectrum $\sigma(A_t)$ is Fell continuous,

(iii) the spectral edges are continuous.
\end{theo}
%%%%%%%%%%%%%%%%

\noindent Such results can be made more quantitative if a metric is introduced to describe the topology of $T$. Let $(X,d_X)$ and $(Y,d_Y)$ be complete metric spaces. For $\alpha>0$, a function $f:X\to Y$ is called $\alpha$-H\"older whenever
%%%%%%%%%%%%%%%%
$$\Hol^\alpha(f)=
  \sup_{x\neq y}\frac{d_Y(f(x),f(y))}{d_X(x,y)^\alpha}< \infty\,.
$$
%%%%%%%%%%%%%%%%
\noindent A family $\Ff$ of $\alpha$-H\"older functions from $X\to Y$ is called {\em uniformly} $\alpha$-H\"older if $\Hol^\alpha(\Ff)=\sup_{f\in\Ff}\Hol^\alpha(f) <\infty$. Given $M>0$ let $\ps_2(M)$ be the set of polynomials of the form $p(z)=p_0+p_1 z+p_2 z^2$ with $p_i\in \RM$ and $\|p\|_1=|p_0|+|p_1|+|p_2|\leq M$. 

%%%%%%%%%%%%%%%%
\begin{defini}
\label{ccs15.def-p2H}
Let [F1-F3] hold with $(T,d)$ a complete metric space. The field $A$ of self-adjoint bounded operators is called p2-$\alpha$-H\"older continuous  whenever, for all $M>0$, the maps $\Phi_p: t\in T\mapsto \|p(A_t)\|,\; p\in \ps_2(M)$, are uniformly $\alpha$-H\"older. 
\end{defini}
%%%%%%%%%%%%%%%%

%%%%%%%%%%%%%%%%
\begin{rem}
\label{ccs15.rem-Lip}
{\em Lipschitz continuity corresponds to $\alpha=1$.
}%%
\hfill $\Box$
\end{rem}
%%%%%%%%%%%%%%%%

\noindent Theorem~\ref{ccs15.th-main1} indicates that a quantitative control on the norms provides an estimate on the behavior of the spectrum and conversely. This result can be expressed as follows.

%%%%%%%%%%%%%%%%
\begin{theo}
\label{ccs15.th-main2}
Let [F1-F3] hold with $(T,d)$ a complete metric space. Let $A=(A_t)_{t\in T}$ be a family of self-adjoint, bounded operators such that $\sup_{t\in T}\|A_t\|<\infty$.

\begin{itemize}
\item[(i)] If $A$ is a p2-$\alpha$-H\"older continuous field then the spectrum $\sigma(A_t)$ is $\alpha/2$-H\"older continuous with respect to the Hausdorff metric. 

\item[(ii)] If the spectrum $\sigma(A_t)$ is $\alpha$-H\"older continuous with respect to the Hausdorff metric then $A$ is a p2-$\alpha$-H\"older continuous field.
\end{itemize}
\end{theo}
%%%%%%%%%%%%%%%%

\noindent The change in the H\"older exponent in statement (i) is coming from the possible closing of a gap at some point $t_0\in T$. An example can be found in the spectrum of the Almost Mathieu model (eq.~(\ref{ccs15.eq-am})) at the spectral value $c=0$ whenever $t=1/2$ \cite{BS82}. Indeed, near the value $t=1/2$ many gaps have width $O(\sqrt{|t-1/2|})$ \cite{RB90}, while the field $(H_t)_{t\in T}$ is p2-Lipschitz continuous \cite{Be94}. It follows that the proof of this theorem requires a precise definition of gap and gap edges, as well as the concept of gap closing. While the intuitive definition of the gaps and gap edges are actually correct, the definition of gap closing is more subtle that it looks at first and it requires some technicalities. This is explained in Section~\ref{ccs15.ssect-gapedge}. The previous result is made more precise as follows.

%%%%%%%%%%%%%%%%
\begin{theo}
\label{ccs15.th-main3}
Let [F1-F3] hold with $(T,d)$ a complete metric space. Let $A=(A_t)_{t\in T}$ be a family of self-adjoint, bounded operators such that $\sup_{t\in T}\|A_t\|<\infty$. If $A$ is p2-$\alpha$-H\"older continuous then the edges of a spectral gap of $A$ are $\alpha$-H\"older continuous at $t_0$ if this gap does not close at $t_0$.
\end{theo}
%%%%%%%%%%%%%%%%

\noindent It may happen that at some point $t_0\in T$ a spectral gap of $A_t$ closes at the spectral value $c$ called a gap tip. The precise definition of a closing gap and of its tip is given in Definition~\ref{ccs15.def-closedgap}. Such a gap tip will be called {\em isolated} if the distance of $c$ from any gap in the spectrum of $A_{t_0}$ is positive. 

%%%%%%%%%%%%%%%%
\begin{theo}
\label{ccs15.th-main4}
Let [F1-F3] hold with $(T,d)$ a complete metric space. Let $A=(A_t)_{t\in T}$ be a p2-$\alpha$-H\"older continuous field of bounded, self-adjoint operators such that $\sup_{t\in T}\|A_t\|<\infty$. If $c$ is an isolated  gap tip of $A_{t_0}$ then the widths of the gaps of $A_t$ closing on $c$ at $t_0$ are $\alpha/2$-H\"older continuous at $t_0$.
\end{theo}
%%%%%%%%%%%%%%%%

\noindent If the gap tip is not isolated, it is possible to construct examples for which the width of the closing gap can vanish as slowly as one wishes (see Section~\ref{ccs15.ssect-main4}, Example~\ref{ccs15.exam-slowclosing}).

\vspace{.1cm}

\noindent The previous results apply as well for a field $U=(U_t)_{t\in T}$ of unitary operators.

%%%%%%%%%%%%%%%%
\begin{defini}
\label{ccs15.def-unitary}
Let $U=(U_t)_{t\in T}$ be a field of unitary operator. It will be called p-continuous at $t_0$ whenever $\Phi_p:t\in T\mapsto \|p(U_t)\|$ is continuous at $t_0$ for all polynomial of the form $p(z) = 1+e^{\imath\theta}z$ where $\theta\in\TM$. It will be called p-$\alpha$-H\"older, if in addition, this maps are $\alpha$-H\"older uniformly with respect to $\theta$.
\end{defini}
%%%%%%%%%%%%%%%%

%%%%%%%%%%%%%%%%
\begin{theo}
\label{ccs15.th-main5}
Let $U=(U_t)_{t\in T}$ be a field of unitary operator. Then Theorem~\ref{ccs15.th-main1} and Theorem~\ref{ccs15.th-main2} apply, provided p2-continuity is replaced by p-continuity. 
\end{theo}
%%%%%%%%%%%%%%%%

\noindent The proof of this results follows the same lines as for the self-adjoint case. It will be left to the reader.

\vspace{.3cm}

\noindent {\bf Acknowledgments: }S.~B. would like to thank Daniel Lenz for constant support and fruitful discussions. He wants to thank J. Fillman for rising up the case of unitary operators. He also thanks the School of Mathematics at Georgia Institute of Technology for support during his stay in Atlanta in March 2015 to finish this paper. Both authors are thankful to the Erwin Schr\"odinger Institute, Vienna, for support during the Summer of 2014 where parts of this result were obtained. J.~B. thanks Daniel Lenz, the Department of Mathematics and the Research Training Group (1523/2), at the Friedrich-Schiller-University of Jena, Germany for an invitation in May 2015 during which part of this paper was completed. J.~B. also thanks D. Damanik for the Reference \cite{DG11}.

%%%%%%%%%%%%%%%%%%%%%%%%%%%%%%%%%%%%%%%%%%%%%%%%%%%%%%%%%%%%%%%%%%%%
\section{Continuity}
\label{ccs15.sect-pf1}

%%%%%%%%%%%%%%%%%%%%%%%%%%%%%%%%%%%%%%%%%%%%%%%%%%%%%%%%%%%%%%%%%%%%
 \subsection{The Core Argument}
 \label{ccs15.ssect-core}

\noindent The main argument can be phrased as follows. Let $t_0\in T$ and let $(a,b)\subset \RM$ be a gap in the spectrum of $A_0=A_{t_0}$. Namely, $a,b\in\sigma(A_0)$ but $(a,b)\cap \sigma(A_0)=\emptyset$. In order to prove that $b$ becomes a continuous function of $t$ near $t_0$, let $c$ be chosen close to $b$ in the gap, that is $(a+b)/2<c<b$. Then $(b-c)^2$ is the lowest point in the spectrum of $(A_0-c)^2$. If $m>0$ is large enough so that $(A_t-c)^2 < m^2$ for $t$ near $t_0$, it follows that the norm of $m^2-(A_0-c)^2$ is precisely given by $m^2-(b-c)^2$. This norm is continuous in $t$ near $t_0$ by assumption. Following the argument above backward, it follows that, for any $\epsilon$ satisfying $b-c> \epsilon >0$, there is a neighborhood $U$ of $t_0$ such that $(c-\epsilon,c+\epsilon)\cap \sigma(A_t)=\emptyset$ for $t\in U$ and there is $b_t\in\sigma(A_t)$ such that $|b_t-b|<\epsilon$. Such an idea can be expressed more precisely as follows.

%%%%%%%%%%%%%%%%
\begin{lemma}
\label{ccs15.lem-bpA}
Let $A:\Hh\to\Hh$ be a self-adjoint, linear, bounded operator on a Hilbert space $\Hh$. Let $p$ denotes the polynomial defined by $p(z)=m^2-z^2$ with $m > \|A\|$. Then, if $B_r(x)$ denotes the open ball centered at $x$ with radius $r$,

\begin{itemize}
 \item[(i)] the inequality $\|p(A)\|\leq m^2-r^2$, with $m>r$, holds if and only if $B_r(0)\cap\sigma(A)=\emptyset$, and

 \item[(ii)] the inequality $\|p(A)\|> m^2-r^2$, with $m>r$ holds if and only if $B_r(0)\cap\sigma(A)\neq\emptyset$.
\end{itemize}
\end{lemma}
%%%%%%%%%%%%%%%%

\noindent  {\bf Proof: }The two statements are equivalent and only (i) will be proved. Since $A$ is self-adjoint, its spectrum is contained in the real line. Moreover, if $q$ is any polynomial, $\sigma(q(A))=q(\sigma(A))$. It follows that $p(\lambda)\geq m^2-\|A\|^2>0$ for $\lambda\in \sigma(A)$.

\vspace{.1cm}

\noindent Now $B_r(0)\cap\sigma(A)=\emptyset$ if and only if all $\lambda \in \sigma(A)$ satisfy $|\lambda|\geq r$, or, equivalently,  $m^2-\lambda^2\leq m^2-r^2$. This is equivalent to $\|p(A)\|=\sup_{\lambda\in\sigma(A)}\{m^2-\lambda^2\}\leq m^2-r^2$.
\hfill $\Box$

\vspace{.2cm}

\noindent For a unitary operator the analogous result is expressed as follows (the proof is left to the reader).

%%%%%%%%%%%%%%%%
\begin{lemma}
\label{ccs15.lem-bpU}
Let $U:\Hh\to\Hh$ be a unitary, linear operator on a Hilbert space $\Hh$ and $\theta\in\TM$. Then for any $r<2$, the inequality $\|1+e^{-\imath\theta}U\|\leq \sqrt{4-r^2}$ holds if and only if $B_r(e^{\imath\theta})\cap\sigma(U)=\emptyset$.
\end{lemma}
%%%%%%%%%%%%%%%%

%%%%%%%%%%%%%%%%%%%%%%%%%%%%%%%%%%%%%%%%%%%%%%%%%%%%%%%%%%%%%%%%%%%%
 \subsection{Topologies on the Set of Closed Subsets}
 \label{ccs15.ssect-viet}

\noindent Vietoris \cite{Vi22} introduced a topology on the space of closed subsets of a topological space $X$  which is described below (see also \cite{CV77}). Whenever $X$ is a complete metric space, the topology defined by the Hausdorff metric \cite{Ku68,CV77} induces the Vietoris topology. If, furthermore, $X$ is compact the Vietoris topology coincides with the Fell topology \cite{Fe62}.

\vspace{.1cm}

\noindent Here $\cs(X)$ denotes the set of all closed subsets of $X$. For $K\subseteq X$ closed and for $\Ff$ a finite family of open subsets of $X$, let

%%%%%%%%%%%%%%%%
$$\Uu(K,\Ff):=
   \{F\in\cs(X)\,;\, K\cap F=\emptyset,\; 
      F\cap O\neq \emptyset \text{ for all } O\in\Ff \}\,.
\hspace{1cm}
\mbox{\rm\bf (Hit-and-Miss)}
$$
%%%%%%%%%%%%%%%%

\noindent Then $\ws\subset \cs(X)$ is Vietoris open if it is a union of $\Uu(K,\Ff)$'s. The Fell topology is defined in the same way provided the set $K$ is required to be compact instead of closed. The spectrum of a self-adjoint, bounded operator is a compact subset of $\RM$. For compact sets the Vietoris continuity of the spectra is equivalent to the Fell continuity if restricted to a compact set. In particular, $(F_t)_{t\in T}$ a family of compact sets is Vietoris continuous at $t_0$ if and only if there is a compact $F\subset\RM$ and an open neighborhood $U_0$ of $t_0$ such that $F_t\subset F$ for all $t\in U_0$ and $(F_t)_{t\in U_0}$ is Fell continuous in $\cs(F)$. The fact that the spectra stays in a compact set $F$ follows from the p2-continuity as the operator norm is continuous. Let $\overline{\RM}$ be the completion of $\RM$ with two points $\pm\infty$ at infinity. Then the following is an example of an application for the Vietoris topology

%%%%%%%%%%%%%%%%
\begin{proposi}
\label{ccs15.prop-maxViet}
The maps $F\in\cs(\RM)\mapsto \sup{F}\in \overline{\RM}$ and $F\in\cs(\RM)\mapsto \inf{F}\in \overline{\RM}$ are Vietoris continuous.
\end{proposi}
%%%%%%%%%%%%%%%%

\noindent  {\bf Proof: }Since multiplication by $-1$ is Vietoris continuous, the Vietoris continuity of $\sup\{\cdot\}$ leads to Vietoris continuity of $\inf\{\cdot\}$. It is therefore sufficient to prove the Vietoris continuity of $\sup\{\cdot\}$. Let $F_0\in \cs(\RM)$ and let $\lambda:=\sup{F_0}$ be finite. Since $F_0$ is closed, $\lambda\in F_0$. For $\epsilon >0$ let $K_\epsilon$ and $O_\epsilon$ be defined as the closed subset $K_\epsilon:=[\lambda+\epsilon,\infty)$ and the open subset $O_\epsilon:=(\lambda-\epsilon,\infty)$. By construction, $F_0\in \Uu(K_\epsilon, \{O_\epsilon\})$. Now, if $F\in \Uu(K_\epsilon, \{O_\epsilon\})$, it follows that $F\cap K_\epsilon=\emptyset$ which implies $\sup{F}<\lambda +\epsilon$. Moreover, since $F\cap O_\epsilon\neq \emptyset$, it follows that $\lambda-\epsilon <\sup{F}$. Consequently, $F\in \Uu(K_\epsilon, \{O_\epsilon\})$ implies $|\sup {F}-\sup {F_0}|< \epsilon$, proving the continuity at $F_0$. If now $\sup{F_0}=\infty$, the same argument works with $K=\emptyset$ and $O_c:=(c,\infty)$ for any $c>0$.
\hfill $\Box$

\vspace{.2cm}

\noindent As a reminder, a function $f:X\to Y$ between two topological spaces is called {\em closed} if the image under $f$ of any closed subset of $X$ by $f$ is closed in $Y$. 

%%%%%%%%%%%%%%%%
\begin{proposi}{(See \cite{Fi13})}
\label{ccs15.prop-cfc}
Let $X,Y$ be locally compact spaces and let $f:X\to Y$ be a continuous closed function. Then, the map $\hf:F\in\cs(X)\mapsto f(F)\in\cs(Y)$ is Vietoris continuous.
\end{proposi}
%%%%%%%%%%%%%%%%

\noindent {\bf Proof: } Since $f$ is closed, the map $\hf$ is well defined. Let now $K\subset Y$ be closed and let $\Ff$ be a finite family of open subsets of $Y$. Hence $\Uu(K,\Ff)$ is an open set in $\cs(Y)$. Let $\ws=\Uu(f^{-1}(K),f^{-1}(\Ff))$. Since $f$ is continuous, $f^{-1}(O)$ is open for any $O\in \Ff$ and $f^{-1}(K)$ is closed. Hence $\ws$ is an open set in $\cs(X)$. Then $F\in \ws$ means $F\cap f^{-1}(K)=\emptyset$ while $F\cap f^{-1}(O)\neq \emptyset$ for all $O\in\Ff$. This implies that $f(F)\cap K=\emptyset$. For indeed otherwise, there would be an $x\in F$ with $f(x)\in K$, implying that $x\in f^{-1}(K)$, a contradiction. In addition, if $O\in\Ff$ and if $x\in F\cap f^{-1}(O)\neq \emptyset$ then $f(x)\in f(F)\cap O$ follows showing that $f(F)\cap O\neq \emptyset$. Hence $\hf(F)\in \Uu(K,\Ff)$, which implies $\ws\subset {\hf}^{-1}(\Uu(K,\Ff))$, namely, $\hf$ is continuous.
\hfill $\Box$

\vspace{.2cm}

\noindent If $(X,d)$ is a complete metric space it turns out that the Vietoris topology on $\cs(X)$ coincides with the topology induced by the {\em Hausdorff metric} \cite{Ku68,Mu75,CV77,Ba88}. The Hausdorff metric $d_H$ is defined as follows: if $x\in X$ and $A\subset X$, then $\dist(x,A)=\inf\{d(x,y)\,;\, y\in A\}$. Given two subsets $A,B\subset X$, $\delta(A,B)$ is defined by $\delta(A,B)=\sup_{x\in A}\dist(x,B)$. Then the Hausdorff distance of the two sets $A,B$ is simply $d_H(A,B)=\max\{\delta(A,B),\delta(B,A)\}$. In general $d_H$ is only a pseudo-metric on the set of all subsets of $X$. However if restricted to the set $\cs(X)$ it becomes a metric defining the Vietoris topology.

%%%%%%%%%%%%%%%%%%%%%%%%%%%%%%%%%%%%%%%%%%%%%%%%%%%%%%%%%%%%%%%%%%%%
 \subsection{Proof of Theorem~\ref{ccs15.th-main1}}
 \label{ccs15.ssect-main1}

\noindent Let $A=(A_t)_{t\in T}$ be a field of self-adjoint, bounded operators satsfying [F1-F3]. %(see Definition~\ref{ccs15.def-p2}). Let $t_0\in T$. By continuity of the norm $\|A_t\|$ with respect to $t$, for any $m>\|A_{t_0}\|$  there is an open neighborhood $U_0$ of $t_0$ such that if $t\in U_0$ then $\|A_t\|< m$.

%%%%%%%%%%%%%%%%
\begin{lemma}
\label{ccs15.lem-VietSp}
If $A=(A_t)_{t\in T}$ is a p2-continuous field of self-adjoint, bounded operators, then the map $t\in T\mapsto \sigma(A_t)\subset \RM$ is continuous in the Vietoris topology.
\end{lemma}
%%%%%%%%%%%%%%%%

\noindent  {\bf Proof: } Let $t_0\in T$. By continuity of the norm $\|A_t\|$ with respect to $t$, there is an open neighborhood $U_0$ of $t_0$ such that if $t\in U_0$ then $\|A_t\|< \|A_{t_0}\|+1/2$. Set $c:=\|A_{t_0}\|+1$ so that, for $t\in U_0$, $\sigma(A_t)\subset [-c,+c]$. Let $K\subset \RM$ be closed and $\Ff$ be a finite set of open subsets of $\RM$ chosen so that $\sigma(A_{t_0})\in \Uu(K,\Ff)$. Let $K_c:=K\cap[-c,c]$, which is compact. Since $K_c$ and $\sigma(A_{t_0})$ are closed and since $K_c\cap\sigma(A_{t_0})=\emptyset$, for any given $x\in K_c$, there exists an $r(x)>0$ so that $B_{r(x)}(x)\cap \sigma(A_{t_0})=\emptyset$. The family of (smaller) open balls $\{B_{r(x)/2}(x)\,;\, x\in K_c\}$ covers $K_c$. By compactness of $K_c$, there is a finite set $\{x_1,\cdots, x_l\}\subset K_c$ such that

%%%%%%%%%%%%%%%%
$$K_c\subset \bigcup_{k=1}^l B_{r_k/2}(x_k)\,,
   \hspace{2cm}
    B_{r_k}(x_k)\cap\sigma(A_{t_0})=\emptyset\,,
     \hspace{1cm}
      r_k:=r(x_k),\;\,  k\in\NM\,.
$$
%%%%%%%%%%%%%%%%

\noindent Let now $m$ be chosen so that $2\sup_{t\in U_0}\|A_{t}\|+ \sup_{x\in K_c}|x| <m$. Using Lemma~\ref{ccs15.lem-bpA}, the condition $B_{r_k}(x_k)\cap\sigma(A_{t_0})=\emptyset$ is equivalent to $\|m^2-(A_{t_0}-x_k)^2\|< m^2-r_k^2$. The p2-continuity implies that there is an open neighborhood $U_k$ of $t_0$ such that for $t\in U_k$ then $\|m^2-(A_{t}-x_k)^2\|< m^2-r_k^2/4$. The set $U=\bigcap_{k=1}^l U_k\cap U_0$ is open and contains $t_0$. It follows from the previous bounds and from Lemma~\ref{ccs15.lem-bpA}, that, for $t\in U$, $B_{r_k/2}(x_k)\cap\sigma(A_{t})=\emptyset$ for all $k\in\{1,2,\cdots,l\}$ and $\sigma(A_{t})\cap K_c=\emptyset$ implying $K\cap\sigma(A_{t})=\emptyset$.

\vspace{.1cm}

\noindent Similarly, let $O\in\Ff$. Since $O\cap \sigma(A_{t_0})\neq \emptyset$, it follows that for any $x\in O\cap \sigma(A_{t_0})$ there is an $r(x)=r$ such that $B_r(x)\subset O$ (since $O$ is open in $\RM$). Since $x\in \sigma(A_{t_0})$, then $|x|\leq \|A_{t_0}\|$ so that $\|A_{t_0}-x\|\leq 2\|A_{t_0}\|<m$. Consequently, since $B_{r/2}(x)\subset O$ and $x\in\sigma(A_{t_0})$, Lemma~\ref{ccs15.lem-bpA} implies that $\|m^2-(A_{t_0}-x)^2\|>m^2-r^2/4$. Using the p2-continuity, there is an open neighborhood $V_O$ of $t_0$ in $T$ such that for $t\in V_O$ the inequality $\|m^2-(A_{t}-x)^2\|>m^2-r^2$ holds. This leads to $O\cap \sigma(A_t)\supset B_{r}(x)\cap \sigma(A_t)\neq \emptyset$ for $t\in V_O$. As the family $\Ff$ is finite, the intersection $V=\bigcap_{O\in\Ff} V_O$ is open and contains $t_0$. Using the first part of the proof, $V\cap U$ is  an open neighborhood of $t_0$ as well. Consequently, for $t\in U\cap V$ the spectrum of $A_t$ satisfies $\sigma(A_t)\in\Uu(K,\Ff)$.
\hfill $\Box$

%%%%%%%%%%%%%%%%
\begin{lemma}
\label{ccs15.lem-fsp}
Let $A=(A_t)_{t\in T}$ be a field of self-adjoint, bounded operators satisfying the assumptions [F1-F3]. Let $f:\RM\to \CM$ be continuous and closed. If the spectrum map $t\in T\mapsto \sigma(A_t)$ is Vietoris continuous, then the norm-map $t\in T\mapsto \|f(A_t)\|\in \RM_+$ is continuous.
\end{lemma}
%%%%%%%%%%%%%%%%

%%%%%%%%%%%%%%%%
\begin{rem}
\label{ccs15.rem-sidp2}
{\em Thanks to Lemma~\ref{ccs15.lem-fsp}, $A=(A_t)_{t\in T}$ is p2-continuous if $t\in T\mapsto \sigma(A_t)$ is Vietoris continuous since any polynomial $p:\RM\to\CM$ is continuous and closed.
}%%
\hfill $\Box$
\end{rem}
%%%%%%%%%%%%%%%%

\noindent  {\bf Proof: }As $A_t$ is self-adjoint and bounded, $\sigma(A_t)$ is compact in $\RM$ and $f(\sigma(A_t))=\sigma(f(A_t))\subset \CM$. Moreover, $f(A_t)$ is a normal operator, so that the Spectral Theorem applies. The norm map can be seen as the composition of the following continuous maps
%%%%%%%%%%%%%%%%
$$t\mapsto \sigma(A_t)
    \stackrel{\hf}{\mapsto} \sigma(f(A_t))
     \stackrel{|\cdot|}{\mapsto} |\sigma(f(A_t))|
      \stackrel{\sup}{\mapsto} \|f(A_t)\|\,.
$$
%%%%%%%%%%%%%%%%
\noindent In this formula, the map $|\cdot |$ is nothing but $\hg$ if $g:z\in \CM\mapsto |z|\in [0,\infty)$. The leftmost map is continuous by assumption. Thanks to Proposition~\ref{ccs15.prop-cfc}, the second and third maps on the left (with $X,Y=\RM \;\mbox{\rm or}\;\CM$) are continuous. At last Proposition~\ref{ccs15.prop-maxViet} implies that the rightmost map is continuous.
\hfill $\Box$

\newpage
%%%%%%%%%%%%%%%%%%%%%%%%%%%%%%%%%%%%%%%%%%%%%%%%%%%%%%%%%%%%%%%%%%%%
 \subsection{Gaps and Edge Continuity}
 \label{ccs15.ssect-gapedge}

%%%%%%%%%%%%%%%%
\begin{defini}
\label{ccs15.def-gap}
Let $F\subset\RM$ be closed. Any connected component of its complement is an open interval. An edge of $F$ is a boundary of such intervals. Such an interval is called a gap of $F$ if it is bounded.
\end{defini}
%%%%%%%%%%%%%%%%

\noindent It follows that a gap is an open interval $(a,b)$ with $a,b\in F$, $-\infty < a < b <+\infty$ and $(a,b)\cap F=\emptyset$. Then $x\in\RM$ is called a {\em gap edge} of $F$ if there exists a gap $(a,b)$ of $F$ such that either $a=x$ or $b=x$. In addition to the gap edges, the other edges of $F$ are $\sup{F}$ and $\inf{F}$, which might be $\pm \infty$. 

%%%%%%%%%%%%%%%%
\begin{defini}
\label{ccs15.def-closedgap}
Let $T$ be a topological space and let $(F_t)_{t\in T}$ be a family of closed subsets of $\RM$. Given $t_0\in T$, a gap tip $c$ of $F_{t_0}$ is a point $c\in F_{t_0}$ such that there is a non-empty set $U\subset T$ with $t_0\in\overline{U}\setminus U$, and there are two  functions $a:t\in U\mapsto a_t\in \RM$ and $b:t\in U\mapsto b_t\in \RM$ such that $a_t< b_t$, $\lim_{t\to t_0}a_t=\lim_{t\to t_0}b_t=c$ and the interval $(a_t,b_t)$ is a gap of $F_t$.

\noindent A gap tip $c$ of $F_{t_0}$ is called isolated if there is $\delta>0$ such that $(c-\delta,c+\delta) \subset F_{t_0}$. 
\end{defini}
%%%%%%%%%%%%%%%%

\noindent In the previous definition, the family $(F_t)_{t\in T}$ needs not be Vietoris continuous. In addition $U$ needs not be open either. Gap tips occur in various ways. For example, the Almost Mathieu model (see eq.~(\ref{ccs15.eq-am})) admits an isolated gap tip at $t=1/2$ and at the spectral parameter $c=0$ \cite{BS82}. Similarly, the Kohmoto model, a model similar to the Almost Mathieu one, with the potential $V_t(n)=\cos{2\pi(nt+\theta)}$ replaced by $V_t(n)=\chi_{[0,t)}(nt+\theta)$ (where $\chi_I$ denotes the characteristic function of the interval $I\subset \TM$), admits gap tips at any rational values of $t$ \cite{OK85,BIT91}, corresponding to $c$ being an isolated eigenvalue, namely $c$ is the boundary of a left gap $(a,c)$ and of a right one $(c,b)$

\vspace{.1cm}

\noindent Another way to describe the continuity of the spectrum consists in using the continuity of its edges. This concept is a bit delicate, in particular because of the possibility of closing gaps defined in (G3). Here is a definition,

%%%%%%%%%%%%%%%%
\begin{defini}
\label{ccs15.def-contgapedges}
Let $T$ be a topological space and $(F_t)_{t\in T}$ be a family of compact subsets of $\RM$. Then the gap edges of $(F_t)_{t\in T}$ are called {\em continuous at $t_0\in T$} if all of the three following assertions hold.

\begin{itemize}
\item[(G1)] The supremum $t\in T\mapsto\sup{F_t}$ and the infimum $t\in T\mapsto\inf{F_t}$ are continuous at $t=t_0$.

\item[(G2)] Let $(a,b)$ be a gap of $F_{t_0}$ and $\epsilon>0$ be arbitrary. Then there exists an open neighborhood $U:=U(\epsilon,a,b)$ of $t_0$ such that for each $t\in U$ there is a gap $(a_t,b_t)$ of $F_t$ satisfying $|a-a_t|<\epsilon$ and $|b-b_t|<\epsilon$.

\item[(G3)] Let $U\subset T$ be a subset such that $t_0\in \overline{U}\setminus U$. If there are maps $a:U\to \RM$ and $b:U\to \RM$ such that (a) $a_t< b_t$ for $t\in U$, (b) $(a_t,b_t)$ is a gap of $F_t$ for $t\in U$ and (c) $\lim_{t\to t_0}a_t=a$, $\lim_{t\to t_0}b_t=b$, then $(a,b)$ is a gap of $F_{t_0}$ which is closed if $a=b$.
\end{itemize}

\noindent If the gap edges are continuous at any $t_0\in T$ the gap edges are called {\em continuous}.
\end{defini}
%%%%%%%%%%%%%%%%

\noindent It is worth noticing that the set $U$ in (G3) need not be open. This definition leads to the following lemma

%%%%%%%%%%%%%%%%
\begin{lemma}
\label{ccs15.lem-gedges}
Let $T$ be a topological space and $(F_t)_{t\in T}$ be a family of closed subsets of $\RM$. If $(F_t)_{t\in T}$ is Vietoris continuous then the gap edges of $(F_t)_{t\in T}$ are continuous.
\end{lemma}
%%%%%%%%%%%%%%%%

\noindent {\bf Proof: }(i) Let $t_0\in T$. Condition (G1) follows by Proposition~\ref{ccs15.prop-maxViet}.

\vspace{.1cm}

\noindent (ii) Condition (G2): Let $(a,b)$ be a gap of $F_{t_0}$, and let $\epsilon>0$ satisfy $(b-a)/2>\epsilon>0$. Let the following sets be defined: the two open sets $O_a:=(a-\epsilon,a+\epsilon),\; O_b:=(b-\epsilon,b+\epsilon)$ and the closed set $K:=[a+\epsilon,b-\epsilon]$. Since $(b-a)/2>\epsilon$, the inequality $a+\epsilon<b-\epsilon$ holds, and so $K$ is nonempty. By construction, $\Uu(K,\Ff)$ is a neighborhood of $F_{t_0}$, where $\Ff:=\{O_a,O_b\}$. Since $(F_t)_{t\in T}$ is Vietoris continuous, there exists an open neighborhood $U\subset T$ of $t_0$ such that $F_t\in\Uu(K,\Ff)$ for all $t\in U$. Then for $t\in U$, the maximum $a_t:=\max\{O_a\cap F_t\}$ and the minimum $b_t:=\min\{O_b\cap F_t\}$ exist, as the intersections are nonempty and $F_t\cap K=\emptyset$. In particular, $(a_t,b_t)$ is a gap of $F_t$. Since $a_t\in O_a$ and $b_t\in O_b$, the inequalities $|a-a_t|<\epsilon$ and $|b-b_t|<\epsilon$ hold. 

\vspace{.1cm}

\noindent (iii) Let the hypothesis (G3) in Definition~\ref{ccs15.def-contgapedges} hold. In such a case, $a\in F_{t_0}$. For otherwise, there is an $r>0$ such that $\Uu([a-r,a+r],\{\RM\})$ is a neighborhood of $F_{t_0}$. Since $\lim a_t=a$, there is an open neighborhood $V\subset T$ of $t_0$ such that $t\in V$ implies $|a_t-a|<r/2$, namely $F_t\notin \Uu([a-r,a+r],\{\RM\})$ if $t\in V\cap U$. Since $t_0\in \overline{U}$, it follows that $V\cap U$ is not empty. Hence the family $(F_t)_{t\in V}$ cannot be continuous at $t=t_0$. Similarly, $b\in F_{t_0}$.

\vspace{.1cm}

\noindent If $a=b$, it follows from the Definition~\ref{ccs15.def-closedgap} that $(a,b)$ is a gap tip of $F_{t_0}$. If $a<b$, the interval $(a,b)$ is a gap of $F_{t_0}$. If not, it follows that $(a,b)\cap F_{t_0}\neq \emptyset$. For indeed, let $x\in (a,b)\cap F_{t_0}$ and let $\epsilon >0$ be small enough so that $x-a>2\epsilon$ and $b-x>2\epsilon$. Then $O_\epsilon= (x-\epsilon, x+\epsilon)\subset (a+\epsilon, b-\epsilon)$. Hence $F_{t_0}\in \Uu(\emptyset, \{O_\epsilon\})$. Using the Vietoris continuity, there is an open neighborhood $V_1$ of $t_0$ such that for $t\in V_1$, it holds that $F_t\in \Uu(\emptyset, \{O_\epsilon\})$. On the other hand, there is another open neighborhood $V_2$ of $t_0$ such that $t\in V_2\cap U$ implies $|a_t-a|<\epsilon\;\mbox{\rm and}\; |b_t-b|<\epsilon$. Since $V_1\cap V_2$ is a nonempty open set and since $t_0\in \overline{U}$, it follows that $V_1\cap V_2\cap U\neq\emptyset$. On the other hand, $(a_t,b_t)$ is a gap of $F_t$ so that $F_t\cap O_\epsilon\neq \emptyset$ for $t\in V_1\cap V_2\cap U$, which is a contradiction.
\hfill $\Box$

%%%%%%%%%%%%%%%%
\begin{lemma}
\label{ccs15.lem-gapedges}
Let $T$ be a topological space and $(F_t)_{t\in T}$ be a family of compact subsets of $\RM$. If the gap edges are continuous then $(F_t)_{t\in T}$ is Vietoris continuous.
\end{lemma}
%%%%%%%%%%%%%%%%

\noindent  {\bf Proof: } Let $t_0\in T$ and define $m:=1+\max\{|\max{F_{t_0}}|, |\min{F_{t_0}}| \}$. Let $K$ be a closed subset of $\RM$ and $\Ff=\{O_1,\ldots,O_n\}$ be a finite family of open subsets of $\RM$ such that $F_{t_0}\in \Uu(K,\Ff)$. It suffices to find a neighborhood $W$ of $t_0$ such that $F_t\in\Uu(K,\Ff)$ for $t\in W$.

\vspace{.1cm}

\noindent Thanks to (G1), there is an open set $U_0\ni t_0$ such that both $\inf{F_t}, \sup{F_t}$ belong to $[-m,+m]$, namely $F_t\subset [-m,+m]$, for $t\in U_0$. 

\vspace{.1cm}

\noindent Let $K_m=K\cap [-m,+m]$. Since $K_m\cap F_{t_0}=\emptyset$, for any $x\in K_m$, there is an $r(x)>0$ such that $B_{r(x)}(x)\cap F_{t_0}=\emptyset$. By compactness there are $x_1,\cdots,x_l \in K_m$ with $r(x_j)=r_j$ such that
%%%%%%%%%%%%%%%%
$$K_m\subset
   \bigcup_{j=1}^l B_{r_j/2}(x_j)\,,
    \hspace{1.5cm}
     B_{r_j}(x_j)\cap F_{t_0}=\emptyset\,,
      \hspace{1.5cm}
       1\leq j\leq l\,.
$$
%%%%%%%%%%%%%%%%
\noindent If $1\leq j\leq l$ is fixed, there is an open set $U_j\ni t_0$ contained in $U_0$, such that $F_t\cap B_{r_j/2}(x_j)=\emptyset$ for $t\in U_j$. This follows from (G1) if either $x_j\leq \inf{F_{t_0}}$ or $x_j\geq \sup{F_{t_0}}$. Whenever $\inf{F_{t_0}} < x_j <\sup{F_{t_0}}$, then there is a gap $(a,b)$ of $F_{t_0}$ containing $B_{r_j}(x_j)$. Using the condition (G2), $U_j$ can be chosen so that $F_t$ admits a gap $(a_t,b_t)$ satisfying $|a-a_t|<r_j/2$ and $|b-b_t|<r_j/2$ for $t\in U_j$. Hence $B_{r_j/2}(x_j)\subset (a_t,b_t)$ for $t\in U_j$, so that $F_t\cap B_{r_j/2}(x_j)=\emptyset$. Since $K_m$ is covered by the balls $B_{r_j/2}(x_j)$ for $1\leq j\leq l$ and $F_t\subseteq[-m,+m]$ for $t\in U_0$, the intersection $U_K=\bigcap_{j=1}^l U_j\subset U_0$ is a nonempty open neighborhood of $t_0$ such that $K\cap F_t=\emptyset$ for $t\in U_K$. 

\vspace{.1cm}

\noindent Let now $O\in\Ff$. By assumption, there is $x\in O\cap F_{t_0}$, in particular there is an $r>0$ such that $B_{r}(x)\subset O$. Therefore there is an open set $V_O\ni t_0$ satisfying $F_t\cap O\supseteq F_t\cap B_{r}(x)\neq \emptyset$ for $t\in V_O$. Indeed, if $x$ is either $\inf{F_{t_0}}$ or $\sup{F_{t_0}}$, this follows from (G1). Otherwise, suppose the contrary. Then for each open set $V\ni t_0$, there is $t_V\in V$ such that $F_{t_V}\cap B_{r}(x)= \emptyset$. Let U denote the set of all such $t_V$'s. By construction, $t_0\in \overline{U}\setminus U$.

\vspace{.1cm}

\noindent Without loss of generality, $r>0$ can be chosen so that $B_{r}(x)\subset O\cap (\inf{F_{t_0}},\sup{F_{t_0}})$. Thus there is an $\epsilon >0$ such that $\inf{F_{t_0}}+\epsilon<x-r<x+r< \sup{F_{t_0}}-\epsilon$. Thanks to (G1), there is an open set $V_0\subset U_0$ containing $t_0$ such that $F_t\cap [-m,x-r]$ and $F_t\cap [x+r,+m]$ are nonempty for $t\in V_0$. If
%%%%%%%%%%%%%%%%
$$a_t:=\sup\left\{F_t\cap [-m,x-r]\right\}\,,
   \hspace{2cm}
    b_t:=\inf\left\{F_t\cap [x+r,+m]\right\}\,,
     \hspace{2cm} t\in V_0\,,
$$
%%%%%%%%%%%%%%%%
\noindent then $a_t$ and $b_t$ are elements of $F_t$, while $(a_t,b_t)\cap F_t=\emptyset$ for $t\in V_0\cap U$. All limit points of the $a_t$'s belong to $[-m,x-r]$, and, similarly, all limit points of the $b_t$'s belong to $[x+r,+m]$. Choosing a suitable subset $U'\subset V_0\cap U$, we may assume that the limits
%%%%%%%%%%%%%%%%
$$\lim_{t\in U';t\to t_0}a_t=a\,,
   \hspace{2cm}
    \lim_{t\in U';t\to t_0}b_t=b\,,
$$
%%%%%%%%%%%%%%%%
\noindent exist. In particular $a\leq x-r$ and $x+r\leq b$, so that $(x-r,x+r)\subset (a,b)$. Then, thanks to (G3), with $U$ replaced by $U'$, $(a,b)$ is a gap of $F_{t_0}$. But this is a contradiction since $x\in  F_{t_0}$. Hence the open set $V_O$ exists.

\vspace{.1cm}

\noindent Since $\Ff$ is a finite set, $V_{\Ff}=\bigcap_{O\in\Ff} V_O$ is an open neighborhood of $t_0$ such that $F_t\cap O\neq \emptyset$ for $O\in \Ff$ and $t\in V_{\Ff}$. Hence $W=U_K\cap V_{\Ff}$  is an open neighborhood of $t_0$ such that $F_t\in \Uu(K,\Ff)$ for $t\in W$.
\hfill $\Box$

%%%%%%%%%%%%%%%%%%%%%%%%%%%%%%%%%%%%%%%%%%%%%%%%%%%%%%%%%%%%%%%%%%%%
 \subsection{Proof of Theorem~\ref{ccs15.th-mainR}}
 \label{ccs15.ssect-mainR}

\noindent The main remark is the following: let $z=x+\imath y$ with $y\neq 0$. Then, if $A$ is a self-adjoint linear operator on some Hilbert space
%%%%%%%%%%%%%%%%
$$\left\| 
    \frac{1}{A-z}
  \right\|= 
    \left\| 
    \frac{1}{(A-x)^2+y^2}
  \right\|^{1/2}\,.
$$
%%%%%%%%%%%%%%%%
\noindent If the spectrum of $A$ admits a gap $(a,b)$ and if $(a+b)/2 < x< b$, it follows from the Spectral Theorem that
%%%%%%%%%%%%%%%%
\begin{equation}
\label{ccs15.eq-rescont}
\left\| 
    \frac{1}{A-z}
  \right\|=\left(\frac{1}{(b-x)^2+y^2}\right)^{1/2}\,.
\end{equation}
%%%%%%%%%%%%%%%%

\noindent Similarly if $a<x<(a+b)/2$ the same argument links the norm of the resolvent to $a$.

\vspace{.1cm}

\noindent If now $(A_t)_{t\in T}$ is an R-continuous field of self-adjoint operators, the continuity of the norm of the resolvent and eq.~(\ref{ccs15.eq-rescont}) imply that the gap edges are continuous. A basis of the Fell topology on $\CM$ is given by the sets $\Uu(K,\Ff)$ where $K\subseteq\CM$ is a compact subset and $\Ff$ is a finite family of open subsets of $\CM$. If a family of closed subsets $(F_t)_{t\in\ts}$ has continuous gap edges then $(F_t)_{t\in\ts}$ is Fell continuous. This can be proven following the lines of the proof of Lemma~\ref{ccs15.lem-gapedges}. Conversely, if the spectrum is Fell continuous, an adaption of Lemma~\ref{ccs15.lem-gedges} implies that the gap edges are continuous. The R-continuity of the field $(A_t)_{t\in T}$ then follows from eq.~(\ref{ccs15.eq-rescont}).

%%%%%%%%%%%%%%%%%%%%%%%%%%%%%%%%%%%%%%%%%%%%%%%%%%%%%%%%%%%%%%%%%%%%
\section{H\"older Continuity}
\label{ccs15.sect-holder}

\noindent The previous Section shows that a simple criterion permits to get continuity of the spectrum of a field of self-adjoint operators. However, in many cases, continuity should be supplemented by more quantitative estimates. Namely if approximating a self-adjoint operator by a family, it is often necessary to control the speed of convergence. This can be done, for instance, if the topological space $T$ is equipped with a metric. In some cases the topological space $T$ can be very irregular, like a Cantor set or a fractal set. So a metric is really the minimal structure that can be considered. In a metric space $(X,d)$, which will always be assumed to be complete, a natural set of functions is the space of Lipschitz continuous functions. But sometimes it is convenient to consider H\"older continuous functions instead.

%%%%%%%%%%%%%%%%%%%%%%%%%%%%%%%%%%%%%%%%%%%%%%%%%%%%%%%%%%%%%%%%%%%%
 \subsection{Metrics: a Reminder}
 \label{ccs15.ssect-metr}

\noindent Let $(X,d)$ be a complete metric space. The metric $d$ is called an {\em ultrametric} whenever the triangle inequality is replaced by $d(x,y) \leq \max\{d(x,z),d(z,y)\}$. If $0<\alpha$, let $d^\alpha$ denote the function $(x,y)\in X\times X\mapsto d(x,y)^\alpha$. Then if $\alpha \leq 1$, the inequality $(a+b)^\alpha\leq a^\alpha+ b^\alpha$, whenever $a,b$ are non negative, implies that $d^\alpha$ is a new metric. It defines the same topology, since the balls are the same. If $d$ is an ultrametric, and if $\Phi:[0,\infty)\to [0,\infty)$ is monotone increasing such that $\Phi(0)=0$, then the function $d_\Phi= \Phi\circ d$ is also an ultrametric defining the same family of balls. 

\vspace{.1cm}

\noindent Given $\epsilon >0$ an $\epsilon$-path $\gamma$ joining $x$ to $y$, denoted by $\gamma:x\stackrel{\epsilon}{\to} y$, is an ordered sequence $(x_0=x, x_1, \cdots , x_{n-1}, x_n=y)$ such that $d(x_{k-1},x_k)<\epsilon$ for $1\leq k\leq n$. The length of $\gamma$ is $\ell(\gamma) = \sum_{k=1}^n d(x_{k-1},x_k)$. The topological space $X$ is connected if and only if, given any pair $x,y\in X$ and any $\epsilon >0$, there is an $\epsilon$-path joining them. The metric $d$ is called a {\em length metric} whenever $d(x,y)$ coincides with the minimal length of paths (for any $\epsilon >0$) joining $x$ to $y$ \cite{Gr99}.

\vspace{.1cm}

\noindent Given $(X,d_X)$ and $(Y,d_Y)$ two metric spaces, and given $\alpha >0$, a function $f:X\to Y$ is $\alpha$-H\"older if there is a $C>0$ such that $d_Y(f(x),f(x'))\leq C d_X(x,x')^\alpha$ for any pair of points $x,x'$ in $X$. It follows that $\alpha$-H\"older functions are continuous. The H\"older constant is defined as 
%%%%%%%%%%%%%%%%
$$\Hol^\alpha (f)= 
   \sup_{x\neq x'}
    \frac{d_Y(f(x),f(x'))}{d_X(x,x')^\alpha} \,.
$$
%%%%%%%%%%%%%%%%
\noindent More generally this definition can be more local as follows \cite{Gr99}: if $r>0$ 
%%%%%%%%%%%%%%%%
$$\Hol^\alpha_r(f)(x)= 
   \sup_{x';0<d_X(x,x')<r}
    \frac{d_Y(f(x),f(x'))}{d_X(x,x')^\alpha}\,,
   \hspace{1cm}
    \Hol^\alpha_r(f)=\sup_{x\in X} \Hol^\alpha_r(f)(x)\leq 
     \Hol^\alpha(f)\,.
$$
%%%%%%%%%%%%%%%%
\noindent Clearly this quantity is a non decreasing function of $r$, so that the limit $r\to 0$ exists and is called the {\em $\alpha$-dilation} of $f$ at $x$:
%%%%%%%%%%%%%%%%
$$\dil^\alpha(f)(x)= \lim_{r\downarrow 0} \Hol^\alpha_r(f)(x)\,,
   \hspace{2cm}
    \dil^\alpha(f)=\sup_{x\in X} \dil^\alpha(f)(x)\leq \Hol^\alpha(f)\,.
$$
%%%%%%%%%%%%%%%%
\noindent For $\alpha=1$ H\"older continuous functions are called {\em Lipschitz}, and $\Hol^1$ is denoted by $\Lip$. If $d_X$ is a length metric it follows that $\Lip(f)=\dil(f)$ \cite{Gr99}. If $d_X$ is a length metric and $\alpha >1$ then $\Hol^\alpha(f)<\infty$ if and only if $f$ is a constant function. However, there are spaces for which some non constant functions are $\alpha$-H\"older continuous for some $\alpha >1$. In particular if $X$ is a Cantor set and $d_X$ an ultrametric the characteristic function of any clopen set is $\alpha$-H\"older for any $\alpha >0$.

\vspace{.1cm}

\noindent In this Section, the topology of $T$ is induced by a metric $d$ for which it is complete. The real line or the complex plane, or any of their subsets, will be endowed with the usual Euclidean metric. If $t_0\in T$, the continuity of $\|A_t\|$ implies that there is an open subset $U_0$ containing $t_0$ such that $\sup_{t\in U_0} \|A_t\|<\infty$. Replacing $T$ by $U_0$, if necessary, there is  no loss of generality in assuming that $\sup_{t\in T}\|A_t\|=m <\infty$. Thanks to Definition~\ref{ccs15.def-p2H}, given any $M>0$ the following constant is finite:
%%%%%%%%%%%%%%%%
$$C_M= \sup\{\Hol^\alpha(\Phi_p)\,;\, \|p\|_1\leq M\}.
$$
%%%%%%%%%%%%%%%%

%%%%%%%%%%%%%%%%%%%%%%%%%%%%%%%%%%%%%%%%%%%%%%%%%%%%%%%%%%%%%%%%%%%%
 \subsection{Proof of Theorem~\ref{ccs15.th-main2}}
 \label{ccs15.ssect-main2}

\noindent Equipped with the Hausdorff metric $d_H$ the space of closed subsets $\cs(X)$ is a metric space \cite{CV77}.

%%%%%%%%%%%%%%%%
\begin{proposi}
\label{ccs15.prop-hfh}
Let $(X,d_X)$ and $(Y,d_Y)$ be two metric spaces and $f:X\to Y$ be a $\alpha$-H\"older continuous closed function. Then, the map $\hf:F\in\cs(X)\mapsto f(F)\in\cs(Y)$ is $\alpha$-H\"older continuous.
\end{proposi}
%%%%%%%%%%%%%%%%

{\bf Proof: } As mentioned in Proposition~\ref{ccs15.prop-cfc} the function $\hf$ is well defined since $f$ is continuous and closed. Let $K,L\in\cs(X)$. A short computation leads to
%%%%%%%%%%%%%%%%
$$\dist(f(x),f(L))=
   \inf\{ d_Y(f(x),f(y))\,;\, y\in L\} \leq
    C\,\inf\{ d_X(x,y)^\alpha\,;\, y\in L \}=
     C\, \dist(x,L)^\alpha.
$$
%%%%%%%%%%%%%%%%
\noindent Maximizing over $x$ and exchanging the roles of $K$ and $L$ establishes the result.
\hfill$\Box$

%%%%%%%%%%%%%%%%
\begin{lemma}
\label{ccs15.lem-specH}
Let $A=(A_t)_{t\in T}$ be p2-$\alpha$-H\"older continuous field of self-adjoint, bounded operators such that $m:=\sup_{t\in T}\|A_t\|<\infty$. Then the spectrum $\sigma(A_t)$ is $\alpha/2$-H\"older continuous with H\"older constant less than $\sqrt{C_{4m^2+2}}$.
\end{lemma}
%%%%%%%%%%%%%%%%

\noindent  {\bf Proof: } Let $s,t\in T$. According to the definition of the Hausdorff metric it suffices to show $\dist(\lambda,\sigma(A_t))\leq \,\sqrt{C_{4m^2+2}}\; d(s,t)^{\frac{\alpha}{2}}$ for all $\lambda\in\sigma(A_s)$. Without loss of generality suppose that $\lambda\in \sigma(A_s)\setminus\sigma(A_t)$. Since $\lambda\in\sigma(A_s)$, it follows that $\|4m^2-(A_s-\lambda)^2\|=4m^2$ and $(A_t-\lambda)^2\leq 4m^2$. The polynomial $p(z)=4m^2-(z-\lambda)^2$ has a norm $\|p\|_1=1+2|\lambda|+4m^2-\lambda^2=4m^2+2-(1-|\lambda|)^2 \leq 4m^2+2$. Since $\lambda\not\in\sigma(A_t)$ the norm $\|p(A_t)\|$ is exactly $4m^2-\dist(\lambda,\sigma(A_t))^2$. Consequently,
%%%%%%%%%%%%%%%%
$$\dist(\lambda,\sigma(A_t))^2 = 
   \big| \|p(A_t)\| - \|p(A_s)\| \big|\leq 
    C_{4m^2+2}\, d(s,t)^\alpha\,.
$$
%%%%%%%%%%%%%%%%
\hfill$\Box$

For a self-adjoint, bounded operator $A$, the maximum $\max{|\sigma(A)|}$ is exactly $\|A\|$. Let $A=(A_t)_{t\in T}$ be a field of self-adjoint, bounded operators such that $m:=\sup_{t\in T}\|A_t\|<\infty$. Then $\sigma(A_t)$ is a subset of the compact subset $[-m,m]$ for all $t\in T$.

%%%%%%%%%%%%%%%%
\begin{lemma}
\label{ccs15.lem-specHnorm}
Let $A=(A_t)_{t\in T}$ be a field of self-adjoint, bounded operators such that $m:=\sup_{t\in T}\|A_t\|<\infty$ and the assumptions [F1-F3] are satisfied. If the spectrum $\sigma(A_t)$ is $\alpha$-H\"older continuous with H\"older constant $C$ then $A$ is a p2-$\alpha$-H\"older continuous field.
\end{lemma}
%%%%%%%%%%%%%%%%

\noindent  {\bf Proof: } Let $K:=[-m,m]\subseteq\RM$ be the compact subset such that $\sigma(A_t)\subset K$ for all $t\in T$. Any polynomial $p(z)= p_0 + p_1 z + p_2 z^2 $ restricted to $K$ is Lipschitz continuous with Lipschitz constant $|p_1|+2\,m\,|p_2|$. Recall that the norm-map can be seen as the composition of the following maps
%%%%%%%%%%%%%%%%
$$t\mapsto \sigma(A_t)
    \stackrel{\hp}{\mapsto} \sigma(p(A_t))
     \stackrel{|\cdot|}{\mapsto} |\sigma(p(A_t))|
      \stackrel{\max}{\mapsto} \|p(A_t)\|\,.
$$
%%%%%%%%%%%%%%%%
\noindent The absolute value $|\cdot|:[-m,m]\to [0,m]$ and the $\max:\Kk(\RM)\to\RM$ are both Lipschitz continuous with Lipschitz constant $1$, so that
%%%%%%%%%%%%%%%%
$$\big| \|p(A_t)\|- \|p(A_s)\|\big|
	\leq (|p_1|+2\,m\,|p_2|)\, d_H(\sigma(A_s),\sigma(A_t))
		\leq (|p_1|+2\,m\,|p_2|)\, C\, d(s,t)^\alpha\,,
$$
%%%%%%%%%%%%%%%%
\noindent for all $s,t\in T$ by using Proposition~\ref{ccs15.prop-hfh}. Hence the number $C_M= \sup\{\Hol^\alpha(\Phi_p)\,;\, \|p\|_1\leq M\}$ is finite proving that $A$ is a p2-$\alpha$-H\"older continuous field.
\hfill$\Box$

%%%%%%%%%%%%%%%%%%%%%%%%%%%%%%%%%%%%%%%%%%%%%%%%%%%%%%%%%%%%%%%%%%%%
 \subsection{Proof of Theorem~\ref{ccs15.th-main3}}
 \label{ccs15.ssect-main3}

%%%%%%%%%%%%%%%%
\begin{lemma}
\label{ccs15.lem-normH}
Let $A=(A_t)_{t\in T}$ be p2-$\alpha$-H\"older continuous such that $m:=\sup_{t\in T}\|A_t\|<\infty$. Then the maximum $t\mapsto\max{\sigma(A_t)}$ and the minimum $t\mapsto\min{\sigma(A_t)}$ are $\alpha$-H\"older with H\"older constant less than $C_{1+m}$.
\end{lemma}
%%%%%%%%%%%%%%%%

\noindent  {\bf Proof: }If $a_t=\inf\sigma(A_t)$ and $b_t=\sup \sigma(A_t)$ then the largest of them $|a_t|\vee |b_t|$ coincides with $\|A_t\|$ while the smallest $c=|a_t|\wedge |b_t|$ satisfies $m-c=\|m-A_t\|$ leading immediately to the result.
\hfill $\Box$

%%%%%%%%%%%%%%%%
\begin{lemma}
\label{ccs15.lem-gapH}
Let $A=(A_t)_{t\in T}$ be p2-$\alpha$-H\"older continuous such that $m:=\sup_{t\in T}\|A_t\|<\infty$. Then for an open gap $\gG_{t_0}$, the neighboring gaps $\gG_t$ satisfying $\lim_{t\to t_0}\gG_t=\gG_{t_0}$ are $\alpha$-H\"older with $\alpha$-dilation less than $3C_{(4m^2+2)}/|\gG_{t_0}|$ where $|\gG_{t_0}|$ denotes the gap width. 
\end{lemma}
%%%%%%%%%%%%%%%%

\noindent  {\bf Proof: }Let $t_0\in T$. Let $\gG_{t_0}= (a,b)$ be a gap of $\sigma(A_{t_0})$, so that $-\infty < a< b<+\infty$ and $|\gG_{t_0}|=b-a$. Let $\gG_{t_0}$ be subdivided into six intervals of equal length $r$, and let $c$ be the point located at distance $|\gG_{t_0}|/3$ from $b$. That is $r=(b-a)/6$ and $c=a+4r=b-2r$. Due to Theorem~\ref{ccs15.th-mainR}, there is an open neighborhood $U$ of $t_0$ in $T$, such that for $t\in U$, the spectrum of $A_t$ has a gap $(a_t,b_t)$ with $|a_t-a|<r$ and $|b_t-b|<r$. This implies (i) $b_t-a_t = b_t-b+b-a+a-a_t>6r-2r=4r>0$, (ii) $b_t-c=b_t-b+b-c$ so that $r<b_t-c<3r$, (iii) $c-a_t=c-a+a-a_t>4r-r=3r$. Consequently $c$ is closer to $b_t$ than to $a_t$. Hence if $t\in U$, the infimum of the spectrum of $(A_t-c)^2$ is exactly $(b_t-c)^2$. Since $c$ belongs to the convex hull of $\sigma(A_t)$, it follows that $|c|\leq \|A_t\|\leq m$, so that $(b_t-c)^2\leq (A_t-c)^2\leq 4m^2$. Thus $4m^2-(b_t-c)^2=\|4m^2-(A_t-c)^2\|$ whenever $t\in U$. The polynomial $p(z)=4m^2-(z-c)^2$ has a norm $\|p\|_1=1+2|c|+4m^2-c^2= 2+4m^2-(1-|c|)^2\leq 4m^2+2$. Hence, if $s,t\in U$, this gives $|\|p(A_t)\|-\|p(A_s)\||= |(b_t-c)^2-(b_s-c)^2| =|b_t-b_s||b_t+b_s-2c|$. With the help of the inequality (ii) above, this leads to $|b_t+b_s-2c| >2r$. Consequently 
%%%%%%%%%%%%%%%%
$$s,t\in U
   \hspace{1cm}\Rightarrow\hspace{1cm}
    |b_s-b_t|\leq \frac{3C_{(4m^2+2)}}{(b-a)}\;\;d(s,t)^\alpha\,.
$$
%%%%%%%%%%%%%%%%
\noindent Since, in this argument, the size of $U$ cannot be controlled, this estimate gives only an upper bound on the $\alpha$-dilation at $t_0$ of the gap edge, independently of which gap edges is considered. Changing $c$ into $a+2r=b-4r$ will replace $b_t$ by $a_t$, so that the same argument leads to the same estimate for the dilation of the lower gap edge.
\hfill $\Box$

%%%%%%%%%%%%%%%%%%%%%%%%%%%%%%%%%%%%%%%%%%%%%%%%%%%%%%%%%%%%%%%%%%%%
 \subsection{Proof of Theorem~\ref{ccs15.th-main4}}
 \label{ccs15.ssect-main4}

%%%%%%%%%%%%%%%%
\begin{lemma}
\label{ccs15.lem-clgapA}
Let $A=(A_t)_{t\in T}$ be p2-$\alpha$-H\"older continuous such that $m:=\sup_{t\in T}\|A_t\|<\infty$. Then if $c$ is a isolated spectral gap tip of $A_{t_0}$, then the gaps $(a_t,b_t)$ of $A_t$ closing at $c$ satisfy

%%%%%%%%%%%%%%%%
\begin{equation}
\label{ccs15.eq-clgp}
b_t-a_t\leq 2\;\;\sqrt{C_{(4m^2+2)}}\; \;d(t,t_0)^{\alpha/2}\,.
\end{equation}
%%%%%%%%%%%%%%%% 
\end{lemma}
%%%%%%%%%%%%%%%%

\noindent  {\bf Proof: }In what follows, let $F_t$ be the spectrum of $A_t$. Let $c\in F_{t_0}$ be an isolated gap tip. Namely, according to Definition~\ref{ccs15.def-closedgap}, and since $c$ is isolated, the following assumptions hold

(i) there is $\delta >0$, such that $(c-\delta,c+\delta)\subset F_{t_0}$,

(ii) there is a non-empty set $U$ such that $t_0\in \overline{U}\setminus U$ and for all $t\in U$, the set $F_t$ admits a gap $(a_t,b_t)$ such that $\lim_{t\to t_0}a_t=c=\lim_{t\to t_0}b_t$.

\vspace{.1cm}

\noindent By definition, there is an open neighborhood $V\subset U$ of $t_0$ in $T$ such that, whenever $t\in V\cap U$, then $A_t$ admits a spectral gap $(a_t,b_t)$ with $a_t<b_t$ such that $\max\{|a_t-c|,|b_t-c|\}<\delta$. It follows that $c-\delta< a_t<b_t<c+\delta$. Choosing $\lambda=(a_t+b_t)/2$,  it follows that $\lambda \in F_{t_0}$ and the distance of $\lambda$ to $F_t$ is exactly $(b_t-a_t)/2$. Consequently
%%%%%%%%%%%%%%%%
$$\|4m^2-(A_{t_0}-\lambda)^2\|=4m^2\,,
   \hspace{2cm}
    \|4m^2-(A_{t}-\lambda)^2\|=4m^2-\frac{(b_t-a_t)^2}{4}\,.
$$
%%%%%%%%%%%%%%%%
\noindent Since the polynomial $p(z)= 4m^2-(z-\lambda)^2$ satisfies $\|p\|_1\leq 4m^2+2$ (see the proof of Lemma~\ref{ccs15.lem-specH}), it follows that 
%%%%%%%%%%%%%%%%
$$(b_t-a_t)^2\;\;\leq 4C_{(4m^2+2)} \;\; d(t,t_0)^\alpha\, ,\qquad t\in V.
$$
%%%%%%%%%%%%%%%%
\hfill $\Box$

\vspace{.2cm}

\noindent The gap closing condition is observed in many models. First, in the Almost-Mathieu model $H_t$ (see eq.~\ref{ccs15.eq-am}), the spectrum is a finite union of intervals when $t\in\QM$, since this is a periodic Hamiltonian. So any gap tip is isolated. It has been proved that the gap predicted by the gap labeling Theorem \cite{Be92} at the spectral parameter $c=0$ for $t=p/2q$ is always closed, namely it is actually a gap tip \cite{BS82}. Since it is p2-Lipschitz thanks to \cite{Be94}, the gap width can be at most $1/2$-H\"older. A semi-classical calculation \cite{RB90} validates this prediction. 

\vspace{.1cm}

Another situation where gaps are closing is provided by a small perturbation of the Laplacian on $\ZM$. Namely, for $\lambda \geq 0$ let $H_\lambda$ be defined on $\ell^2(\ZM)$ by 
%%%%%%%%%%%%%%%%
$$H_\lambda\psi(n)=
   \psi(n+1)+\psi(n-1)+\lambda V(n)\psi(n)\,,
$$ 
%%%%%%%%%%%%%%%%
\noindent where $V(n)$ takes on finitely many values. For $\lambda=0$ the spectrum of $H_\lambda$ is the interval $[-2,+2]$. The prediction provided by the {\em Gap Labeling Theorem} \cite{Be92,Be92b}, shows that for certain potentials $V$, gaps may open as $\lambda$ increases from $\lambda =0$. Explicit calculations made on several examples, such as the Fibonacci sequence \cite{SM90,DG11}, Thue-Morse sequence \cite{Be90}, the period doubling sequence \cite{BBG91}, validate that the gap width is at most $O(\lambda^{1/2})$.

\vspace{.2cm}

\begin{figure}[t]
   \centering
\includegraphics[width=12cm]{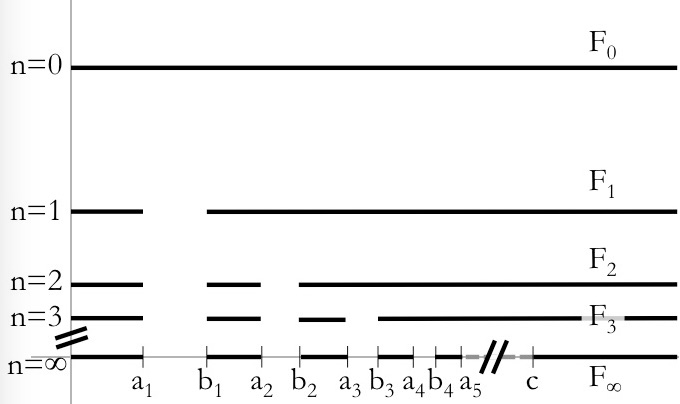}
\caption{Example of a slow closing gap}
\label{ccs15.fig-gapclosing}
\end{figure}

\noindent If the condition that $c$ is isolated is relaxed, the following example (see Fig.~\ref{ccs15.fig-gapclosing}) shows that the closing of gaps cannot be bounded in general. 

%%%%%%%%%%%%%%%%%%
%%%%%%%%%%%%%%%%%%

%%%%%%%%%%%%%%%%
\begin{exam}[{\bf A counter example}]
\label{ccs15.exam-slowclosing}
{\em Let $(a_n,b_n)_{n\in \NM}$ be a double sequence of real numbers such that $0<a_n<b_n<a_{n+1}<c$ and $\sup_{n\in\NM} a_n=c$ (see Fig.~\ref{ccs15.fig-gapclosing}). Let $F_0=[0,m]$ with $m>c$. Then the sequence $(F_n)_{n\in\NM}$ of closed subsets of $F_0$ is defined inductively by $F_{n+1}=F_n\setminus (a_{n+1},b_{n+1})$. Then $F_n$ can be seen as the spectrum of a self-adjoint, bounded operator $A_n$. Clearly $F_n$ converges in the Vietoris and the Fell topologies to $F_\infty=\bigcap_{n\in\NM} F_n$ and $c$ is a gap tip of $F_\infty$. Here $0\leq A_n\leq m$ for all $n$. If $p$ is a polynomial of degree two with real coefficients, then $p(A_n)$ admits $\hat{p}(F_n)$ as its spectrum. Moreover, $p(z)$ can be written as $r+q(z-h)^2$ where $h$ denotes its critical point. Depending on the sign of the coefficients $q,r$, the maximum can be (i) either $p(0)$, (ii) $p(m)$, (iii) $r$, if $h\in F_\infty$, (iv) $p(a_l)$, or $p(b_l)$, whenever $a_l<h<b_l$ and $n\geq l$. Hence $\|p(A_n)\|$ is eventually constant as $n\to \infty$. It follows that if $T=\NM\cup\{\infty\}$ is endowed with any metric such that $d(n,\infty)\to 0$ as $n\to\infty$, the field $A$ is p2-$\alpha$-H\"older for any $\alpha >0$. Let us now consider the case where $|b_{n+1}-a_{n+1}|<|b_n-a_n|$ for all $n$ and let the ultrametric $d$ defined on $\NM\cup\{\infty\}$ as
%%%%%%%%%%%%%%%%
$$d(n,m)= d(m,n)= e^{-\kappa^m}\,,
   \hspace{2cm}
     \mbox{\rm if}\;\; m<n\,, 
$$
%%%%%%%%%%%%%%%%
\noindent and $d(n,n)=0$. Note that $\kappa>1$ is implicitly required to guarantee that $d(n,\infty)\to 0$ if $n\to\infty$. Then $d_H(F_n,F_\infty)=|b_{n+1}-a_{n+1}|/2$. Let the sequences $(a_n),(b_n)$ be chosen such that there is $C>0$ such that
%%%%%%%%%%%%%%%%
$$|b_{n+1}-a_{n+1}|=
   2d_H(F_n,F_\infty)=
    Cd(n,\infty)^{\alpha/2}
$$
%%%%%%%%%%%%%%%%
\noindent It follows that $|b_n-a_n|=Ce^{-(\alpha/2)\kappa^{n-1}}= Cd(n,\infty)^{\alpha/(2\kappa)}$. Hence the gap width is H\"older in $n$ but with an exponent $\alpha/(2\kappa)<\alpha/2$.
}%%
\hfill $\Box$
\end{exam}
%%%%%%%%%%%%%%%% 

\vspace{1cm}
%%%%%%%%%%%%%%%%%%%%%%%%%%%%%%%%%%%%%%%%%%%%%%%%%%%%%%%%%%%%%%%%%%%%%%%%%
%\newpage
%%%%%%%%%%%%%%%%%%%%%%%%%%%%%% July 14 2015 %%%%%%%%%%%%%%%%%%%%%%%%%%%%%% 
%%%%                                                                  %%%%
%%%% Continuity of the spectrum of a field of self-adjoint operators  %%%%
%%%%                                                                  %%%%
%%%%                Siegfried Beckus, Jean Bellissard,                %%%% 
%%%%                   Revised Version #2: 01/21/16                   %%%%
%%%%                                                                  %%%%
%%%%%%%%%%%%%%%%%%%%%%%%%%%%%%%%%%%%%%%%%%%%%%%%%%%%%%%%%%%%%%%%%%%%%%%%%%

%%%%%%%%%%%%%%%%%%%%%%%%%%%%%%%%%%%%%%%%%%%%%%%%%%%%%%%%%%%%%%%%%%%%%%%%%
\end{document}